\documentclass[a4paper]{scrartcl}

\usepackage{amsmath}
\usepackage{amsfonts}
\usepackage{mathtools}
\usepackage{MnSymbol}
\usepackage{epsfig}
\usepackage{graphics}
\usepackage{verbatim} 
\usepackage{siunitx}
\usepackage{subfigure}
\usepackage{booktabs}
\usepackage{authblk}
\usepackage{hyperref}
\hypersetup{
  colorlinks   = true, 
  urlcolor     = red, 
  linkcolor    = blue, 
  citecolor   = blue 
}

\usepackage{algorithmic}
\usepackage{algorithm} 

\newtheorem{definition}{Definition}

\begin{document}

\title{Derivation of higher-order terms in FFT-based numerical homogenization}
\author[1]{ \underline{Felix Dietrich}}
\author[1]{Dennis Merkert}
\author[1]{Bernd Simeon}

\affil[1]{\footnotesize{Technische Universit{\"a}t Kaiserslautern, Paul-Ehrlich-Stra{\ss}e 31\\
D-67663 Kaiserslautern, Germany\\
{\texttt \{\underline{dietrich},dmerkert,simeon\}@mathematik.uni-kl.de}}}

\date {\footnotesize{01.12.2017}}

\maketitle

\begin{abstract}
\textbf{Abstract:} In this paper, we first introduce the reader to the Basic Scheme of Moulinec and Suquet in the setting of quasi-static linear elasticity, which takes advantage of the fast Fourier transform on homogenized microstructures to accelerate otherwise time-consuming computations. By means of an asymptotic expansion, a hierarchy of linear problems is derived, whose solutions are looked at in detail. It is highlighted how these generalized homogenization problems depend on each other. We extend the Basic Scheme to fit this new problem class and give some numerical results for the first two problem orders.
\end{abstract}

\section{Introduction}
\label{dietrich_contrib:sec:1}
Numerical homogenization deals with the 
efficient computation of macroscopic quantities, the so-called
effective properties, by solving microstructural problems in 
representative volume elements (RVEs). Based on the assumption
of a periodic microstructure, efficient FFT-based algorithms \cite{dietrich_contrib:bib:BD2010,dietrich_contrib:bib:MB2012,dietrich_contrib:bib:MS1994} can be applied for this purpose. They have recently been shown to be very competitive, taking advantage of imaging data, i.e. pixels or voxels, as computational mesh.

We study an extension of this methodology that includes higher-order derivatives of the macroscopic quantities, with the aim of 
attaining higher accuracy for the microscopic solutions and the effective properties in this way. This idea has been
introduced by Boutin \cite{dietrich_contrib:bib:B1996}, and we discuss here the algorithmic treatment 
in a unified framework.
For the state-of-the-art in FFT-based numerical homogenization, 
we mention the work on augmented Lagrangians \cite{dietrich_contrib:bib:MMS2001}, the variational scheme based on the Hashin-Shtrikman energy principle \cite{dietrich_contrib:bib:BD2010,dietrich_contrib:bib:BD2012,dietrich_contrib:bib:HS1963}, the polarization scheme \cite{dietrich_contrib:bib:MB2012,dietrich_contrib:bib:MB2013},
and the extension to non-linear problems \cite{dietrich_contrib:bib:CS1997}, elasto-plasticity \cite{dietrich_contrib:bib:SL2014},  elasto-viscoplasticity \cite{dietrich_contrib:bib:EEtAl2013}, and large strains in polycrystals \cite{dietrich_contrib:bib:NEtAl2017}.

The paper is organized as follows.
We start with a compact summary of the so-called Basic Scheme by Moulinec and Suquet  \cite{dietrich_contrib:bib:MS1994,dietrich_contrib:bib:MS1998} and its CG-formulation by Vond\v{r}ejc \cite{dietrich_contrib:bib:VZM2014}.
Then we derive a class of generalized homogenization problems \cite{dietrich_contrib:bib:TMB2012} that include higher-order terms and show how the Basic Scheme can be easily extended for this case.
The paper closes with numerical examples and comparisons.

\section{FFT-based homogenization}
\label{dietrich_contrib:sec:2}
We examine the problem of periodic quasi-static linear elasticity 
in a representative volume element. The $d$-dimensional torus $\mathbb{T}^d\coloneqq\mathbb{R}^d/\mathbb{Z}^d\cong\left[-\frac{1}{2},\frac{1}{2}\right)^d$ is chosen as reference domain to enforce the periodicity. 
For $x\in\mathbb{T}^d$, the microstructure is (in a weak sense) characterized by the strain-displacement equation
\begin{equation}
\epsilon(x)=\frac{1}{2}\left(\nabla u(x)+\Big(\nabla u(x)\Big)^{\text{T}}\right)\enspace,
\label{dietrich_contrib:eq:strain_displacement_eq}
\end{equation}
a constitutive equation in form of Hooke's law
\begin{equation}
\sigma(x)=C(x)\vdotdot\Big(\epsilon(x)+E \Big)\enspace,
\label{dietrich_contrib:eq:Hookes_law}
\end{equation}
with $\vdotdot$ being the double-dot product, and the balance of momentum
\begin{equation}
\nabla\cdot\sigma(x)=0\enspace,
\label{dietrich_contrib:eq:eq_of_motion}
\end{equation}
where we assume that no external forces are applied.
For a given stiffness distribution $C\in L^\infty\left(\mathbb{T}^d\right)^{d\times d\times d\times d}_\text{sym}$ with major and minor symmetries and a prescribed symmetric macroscopic strain tensor $E\in \mathbb{R}^{d\times d}_\text{sym}$, the problem has a unique weak solution for the displacement $u \in H^1(\mathbb{T}^d)^d$, with strain and stress tensors $\epsilon(u), \sigma(u) \in L^2\left(\mathbb{T}^d\right)^{d\times d}_\text{sym}$, such that the mean value $\int_{\mathbb{T}^d}u(x)\,\mathrm{d}x$ is equal to zero; see \cite{dietrich_contrib:bib:KA2013}.

By introducing a regular isotropic reference tensor $C^0\in\mathbb{R}^{d\times d\times d\times d}_\text{sym}$, whose entries are given in the case $d=3$ for the Lam\'{e} parameters $\lambda^0\in\mathbb{R}$ and $\mu^0\in\mathbb{R}\mathbin{\backslash}\{0\}$ as
\begin{equation}
C^0_{ijkl}=\lambda^0\delta_{ij}\delta_{kl}+\mu^0\left(\delta_{ik}\delta_{jl}+\delta_{il}\delta_{jk}\right),\quad i,j,k,l \in \{1,2,3\}, \nonumber
\end{equation}
with $\delta_{ij}$ being the Kronecker-Delta, we can rewrite \eqref{dietrich_contrib:eq:Hookes_law} and \eqref{dietrich_contrib:eq:eq_of_motion} as
\begin{equation}
\nabla\cdot\Big(C^0\vdotdot\epsilon(x)+\tau(x)\Big)=0
\label{dietrich_contrib:eq:auxiliary_problem}
\end{equation}
with the polarization term $\tau(x)\coloneqq\Big(C(x)-C^0\Big)\vdotdot\epsilon(x)+C(x)\vdotdot E$. The solution of \eqref{dietrich_contrib:eq:auxiliary_problem} is given by the periodic Lippmann-Schwinger equation \cite{dietrich_contrib:bib:MS1994,dietrich_contrib:bib:MS1998}
\begin{equation}
\epsilon(x)=-\Big(\Gamma^0\ast\tau\Big)(x)\coloneqq\sum\limits_{\substack{\xi\neq0 \\ \xi\in\left(2\pi\mathbb{Z}\right)^d}}\Big[-\hat{\Gamma}^0(\xi)\vdotdot\hat{\tau}(\xi)\Big]\exp(i\xi\cdot x)\enspace,\label{dietrich_contrib:eq:strain_solution}
\end{equation}
where the Fourier coefficients of the Green strain operator $\Gamma^0$ are given as
\begin{equation}
\hat{\Gamma}^0_{ijkl}(\xi)=\frac{1}{4\mu^0\left\|\xi\right\|^2}(\delta_{ki}\xi_l\xi_j+\delta_{li}\xi_k\xi_j+\delta_{kj}\xi_l\xi_i+\delta_{lj}\xi_k\xi_i)-\frac{\lambda^0+\mu^0}{\mu^0(\lambda^0+2\mu^0)}\frac{\xi_i\xi_j\xi_k\xi_l}{\left\|\xi\right\|^4}\enspace, \nonumber
\end{equation}
for frequencies $\xi\neq0$.

While \eqref{dietrich_contrib:eq:strain_solution} can directly be used for an iterative solution scheme called Basic Scheme, it is preferred to bring all terms containing the strain $\epsilon$ to the left-hand side, resulting in the equation
\begin{equation}
\bigg(\text{Id}+\Gamma^0(x)\ast\Big(C(x)-C^0\Big)\bigg)\vdotdot\epsilon(x)=-\Gamma^0(x)\ast\Big(C(x)\vdotdot E\Big)\enspace.
\label{dietrich_contrib:eq:cg_form}
\end{equation}
The advantage of this formulation is that Krylov subspace methods such as the CG method are directly applicable; see \cite{dietrich_contrib:bib:VZM2014}. The resulting CG-version of the Basic Scheme then reads as follows.
\begin{algorithm}[H]
\caption{Basic Scheme (CG-version)}
\begin{algorithmic}[1]
\STATE \textbf{INIT}:
\STATE \quad $\epsilon_0(x) = -\Gamma^0(x)\ast\Big(C(x)\vdotdot E\Big), \quad \forall x \in\mathbb{T}^d$
\STATE \textbf{ITERATION}:
\STATE \quad $\tau_n(x) = \Big(C(x)-C^0\Big)\vdotdot \epsilon_n(x), \quad \forall x \in\mathbb{T}^d$
\STATE \quad$\hat{\tau}_n = \mathcal{F}(\tau_n)$
\hfill // \, Fourier Transform
\STATE \quad$\hat{\eta}_n(\xi) =  - \hat{\Gamma}^0(\xi)\vdotdot \hat{\tau}_n(\xi), \quad \forall \xi \in \left(2\pi\mathbb{Z}\right)^d\mathbin{\backslash}\left\lbrace0\right\rbrace$
\STATE \quad$\hat{\eta}_{n}(0) = 0$
\STATE \quad$\eta_n = \mathcal{F}^{-1}(\hat{\eta}_n)$
\hfill // \, Inverse Fourier Transform
\STATE \quad$\epsilon_{n+1}(x) = \epsilon_n(x) - \eta_n(x), \quad \forall x \in\mathbb{T}^d$
\STATE \quad\textit{\emph{Check convergence criterion}}
\end{algorithmic}
\label{dietrich_contrib:alg:acc_basic_scheme}
\end{algorithm}
Although different convergence criteria are possible, we will stick to a simple Cauchy criterion of the form $\left\|\epsilon_{n+1}-\epsilon_n\right\|_{L^2}/\left\|\epsilon_0\right\|_{L^2} < \text{TOL}$.

\section{General homogenization problem of order $\alpha$}
\label{dietrich_contrib:sec:3}
Our goal is to extend the problem of quasi-static linear elasticity such that not only macroscopic strains but also macroscopic strain gradients or even higher-order derivatives can be included. This can be achieved by a scale separation as presented in \cite{dietrich_contrib:bib:B1996}. A characteristic length is to be associated with both the macro- and the microscale. The first one will be denoted by $L$, which may be the overall size of the macroscopic sample that is analysed. The latter one is defined by the size of a representative volume element and will be denoted by $\ell$. If the scale ratio $\kappa\coloneqq\ell/L$ is considerably smaller than 1 without being negligible yet, the homogenization framework is applicable. We define the macroscopic variable $Y\coloneqq x/L$ and the microscopic variable $y\coloneqq x/\ell$, which allow for the displacement to be formally written as an asymptotic series expansion
\begin{equation}
u(Y,y) = L\Big(u_0(Y,y) + \kappa u_1(Y,y) + \kappa^2 u_2(Y,y) + \kappa^3 u_3(Y,y) + \ldots\Big)\enspace.
\label{dietrich_contrib:eq:series_displacement}
\end{equation}
In the following, we will usually drop the dependencies on the spatial variables for the sake of better readability. However, if a quantity might only depend on either the microscopic or the macroscopic variable alone, we will denote this explicitly.

By splitting the nabla operator $\nabla = \frac{1}{L}\left( \nabla_Y + \frac{1}{\kappa}\nabla_y\right)$ and by defining symmetric macroscopic and microscopic gradients
\begin{equation}
e_Y(u_i) := \frac{1}{2}\left( \nabla_Y \otimes u_i +  u_i \otimes \nabla_Y\right)
\quad \text{and} \quad
e_y(u_i) := \frac{1}{2}\left( \nabla_y \otimes u_i +  u_i \otimes \nabla_y\right) \nonumber
\end{equation}
for $i=0,1,\ldots,$ accordingly, the series expressions for the strain $\epsilon$ and the stress $\sigma$ can be derived. Inserting \eqref{dietrich_contrib:eq:series_displacement} into the strain-displacement equation \eqref{dietrich_contrib:eq:strain_displacement_eq} gives
\begin{equation}
\epsilon = e_Y (u_0) + \kappa^{-1}e_y (u_0) + \kappa e_Y (u_1) + e_y (u_1) + \ldots
\label{dietrich_contrib:eq:series_strain}
\end{equation}
and after an application of Hooke's law \eqref{dietrich_contrib:eq:Hookes_law} we end up with
\begin{equation}
\sigma = C:e_Y (u_0) + \kappa^{-1}C:e_y (u_0) + \kappa C:e_Y (u_1) + C:e_y (u_1) + \ldots\enspace.
\label{dietrich_contrib:eq:series_stress}
\end{equation}

If we furthermore insert \eqref{dietrich_contrib:eq:series_stress} into the balance of momentum \eqref{dietrich_contrib:eq:eq_of_motion}, we get
\begin{eqnarray}
0 &=& \nabla_Y\cdot\Big[C:e_Y (u_0)\Big] + \kappa^{-1}\nabla_y\cdot\Big[C:e_Y (u_0)\Big] \nonumber\\
&+&\kappa^{-1}\nabla_Y\cdot\Big[C:e_y (u_0)\Big] + \kappa^{-2}\nabla_y\cdot\Big[C:e_y (u_0)\Big] + \ldots\enspace.
\label{dietrich_contrib:eq:series_eq_of_motion}
\end{eqnarray}
Each term $u_i$ appears in four different addends, of which one consists only of purely macroscopic derivatives, one only of purely microscopic derivatives and the remaining two consist of mixed derivatives. We introduce the notation
\begin{eqnarray*}
P^{0}(u_i)&:=&\nabla_Y\cdot\Big[C:e_Y(u_i)\Big]\enspace,\\
P^{-1}(u_i)&:=&\nabla_Y\cdot\Big[C:e_y(u_i)\Big]+\nabla_y\cdot\Big[C:e_Y(u_i)\Big]\enspace,\\
P^{-2}(u_i)&:=&\nabla_y \cdot\Big[C:e_y(u_i)\Big]\enspace.
\end{eqnarray*}
Rearranging the terms of \eqref{dietrich_contrib:eq:series_eq_of_motion} with respect to the exponent of $\kappa$ leads to the expression
\begin{eqnarray*}
0 &=& \kappa^{-2} \Big[P^{-2}(u_0)\Big] + \kappa^{-1} \Big[P^{-2}(u_1) + P^{-1}(u_0) \Big] \nonumber\\
&+& \kappa^0 \Big[P^{-2}(u_2) + P^{-1}(u_1) + P^{0}(u_0)\Big] + \ldots\enspace,
\end{eqnarray*}
where each bracket has to vanish for the left-hand side to be zero. This structure allows us to solve for the term $u_i$ successively in a hierarchical manner.

The first problem $0 = P^{-2}(u_0) = \nabla_y \cdot\left[C:e_y(u_0)\right]$ is trivially solved by a purely macroscopic displacement $u_0(Y,y) = U(Y)$.

The second problem takes the form
\begin{eqnarray*}
0 &=& P^{-2}(u_1) + P^{-1}(u_0) \nonumber\\
&=& \nabla_y \cdot  \Big[ C:e_y (u_1) \Big] + \nabla_Y \cdot  \left[ C:e_y \Big(U(Y)\Big) \right] + \nabla_y \cdot  \left[ C:e_Y \Big(U(Y)\Big) \right] \nonumber\\
&=& \nabla_y \cdot  \Big[ C:e_y (u_1) \Big] + \nabla_y \cdot  \Big[ C:E(Y) \Big] \nonumber\\
&=& \nabla_y \cdot  \left[ C:\Big(e_y (u_1) + E(Y)\Big) \right]\enspace,
\end{eqnarray*}
which coincides with the classical problem presented in Section \ref{dietrich_contrib:sec:2}. Its solution can be computed with Algorithm \ref{dietrich_contrib:alg:acc_basic_scheme}.

All the higher-order problems have essentially the same structure. We restrict ourselves to the second order problem
\begin{eqnarray*}
0 &=& P^{-2}(u_0) + P^{-1}(u_1) + P^0(u_2) \nonumber\\
&=& \nabla_y \cdot  \Big[ C:e_y (u_2) \Big] + \nabla_Y \cdot  \Big[ C:e_y (u_1) \Big] + \nonumber\\
&\phantom{=}& \nabla_y \cdot  \Big[ C:e_Y (u_1) \Big] + \nabla_Y \cdot  \Big[ C:e_Y \Big(U(Y)\Big) \Big]\enspace,
\end{eqnarray*}
but the following idea applies to the remaining higher-order problems as well.
The displacement term $u_1(Y,y)$ depends linearly on the macroscopic strain $E(Y)$. Therefore, we use a separation of variables to make the ansatz $u_1(Y,y)=X_1(y)\vdotdot E(Y)$ with $X_1(y) \in H^1(\mathbb{T}^d)^{d\times d\times d}$ being a third-order tensor depending solely on the microscopic variable, which has to be determined beforehand. The above problem then reads
\begin{eqnarray*}
0 &=& \nabla_y \cdot  \Big[ C:e_y (u_2) \Big] + \nabla_Y \cdot  \Big[ C:e_y \Big(X_1(y)\vdotdot E(Y)\Big) \Big] \nonumber\\
&+& \nabla_y \cdot  \Big[ C:e_Y \Big(X_1(y)\vdotdot E(Y)\Big) \Big] + \nabla_Y \cdot  \Big[ C:e_Y \Big(U(Y)\Big) \Big] \enspace.
\end{eqnarray*}
After rearranging the terms, one can define the polarization
\begin{equation*}
p_2 \coloneqq \frac{1}{2}C\vdotdot\bigg[X_1(y)\vdotdot\nabla E(Y)+\Big(X_1(y)\vdotdot\nabla E(Y)\Big)^\text{T}\bigg]
\end{equation*}
and the body force
\begin{equation*}
g_2 \coloneqq \nabla_Y\cdot\bigg[C\vdotdot\Big[e_y\Big(X_1(y)\vdotdot E(Y)\Big)+e_Y\Big(U(Y)\Big)\Big]\bigg]\enspace.
\end{equation*}
The problem is reduced to the equation
\begin{equation}
0 = \Big[C\vdotdot e_y(u_2) + p_2\Big] + g_2\enspace,
\end{equation}
a general form also taken by the remaining higher-order problems \cite{dietrich_contrib:bib:B1996,dietrich_contrib:bib:TMB2012}. Thus, we define the generalized homogenization problem of order $\alpha$, for $\alpha = 1,2,\ldots,$ as follows; see also \cite{dietrich_contrib:bib:TMB2012}.
\begin{definition}
For $Y\in \Omega$ fixed and $y \in \mathbb{T}^d$, the generalized homogenization problem of order $\alpha$ takes the form
\begin{equation*}
\nabla_y\cdot\left(C\vdotdot\epsilon_\alpha(u_\alpha) + p_\alpha\right) + g_\alpha = 0\enspace,
\end{equation*}
where the polarization
\begin{equation*}
p_\alpha =
\begin{cases}
C\vdotdot E(Y) &, \text{ for } \alpha = 1, \\
\frac{1}{2}C\vdotdot\left[X_{\alpha-1}(y) \cdot \nabla^{\alpha-1}E(Y) + \left( X_{\alpha-1}(y) \cdot \nabla^{\alpha-1}E(Y) \right)^\text{T} \right] &, \text{ for } \alpha \geq 2,
\end{cases}
\end{equation*}
and the body force
\begin{equation*}
g_\alpha =
\begin{cases}
0 &, \text{ for } \alpha = 1, \\
C\vdotdot\Big[ e_y \big(X_{1}(y) \vdotdot E(Y)\big) + e_Y \big(U(Y)\big)\Big] &, \text{ for } \alpha = 2,\\
C\vdotdot\Big[ e_y \big(X_{\alpha-1}(y) \cdot \nabla^{\alpha-1}E(Y)\big) + e_Y \big(X_{\alpha-2}(y) \cdot \nabla^{\alpha-2}E(Y)\big)\Big] &, \text{ for } \alpha \geq 3,
\end{cases}
\end{equation*}
are order-dependent terms. It has a unique weak solution $u_\alpha \in H^1(\mathbb{T}^d)^d$ assuming the displacements have a mean value of zero \cite{dietrich_contrib:bib:TMB2012,dietrich_contrib:bib:VZM2014}.
\end{definition}

The solution can be computed with a slight variation of Algorithm \ref{dietrich_contrib:alg:acc_basic_scheme} presented at the end of Section \ref{dietrich_contrib:sec:2}. The only part that has to be adapted is the initialization in Line 2, whereas the iteration loop remains unchanged. We define the quantity $\theta_\alpha$ which has a closed expression in terms of its Fourier coefficients for non-zero frequencies $\xi$ \cite[eq. (38)]{dietrich_contrib:bib:TMB2012} that reads
\begin{equation}
\hat{\theta}_\alpha(\xi) = \frac{i}{\| \xi \|^4}\left[\left(\xi \otimes \xi\right) \hat{g}_\alpha(\xi) \cdot \xi - \Big( \hat{g}_\alpha(\xi)\otimes \xi + \xi \otimes \hat{g}_\alpha(\xi)\Big) \|\xi\|^2 \right] \enspace.
\end{equation}
The algorithm then reads the same as before with the initialization $\epsilon_0 = -\Gamma^0 \ast \left(p_\alpha + \theta_\alpha\right)$, for all $x \in \mathbb{T}^d$, instead.

\section{Numerical results}
\label{dietrich_contrib:sec:4}
\begin{figure}[bth]
\centering
\includegraphics[scale=1]{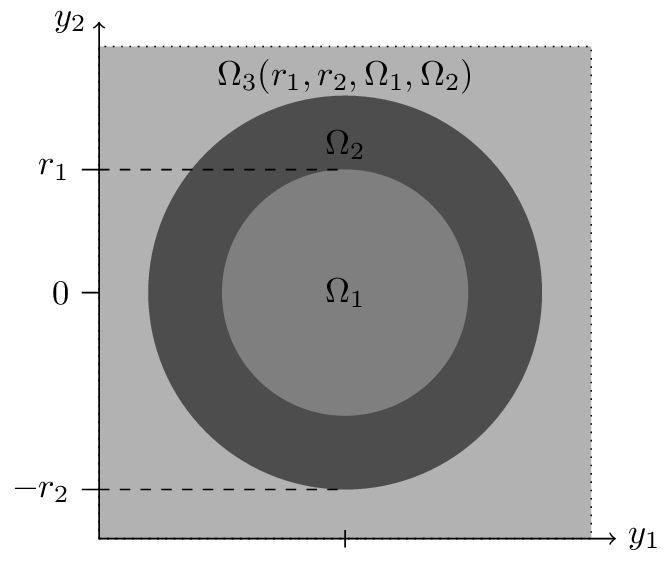}
\caption{Hashin's structure with three distinctive materials $\Omega_i,\, i=1,2,3,$ and radii $0 < r_1 < r_2 \leq 0.5$.}
\label{dietrich_contrib:fig:Hashin}
\end{figure}

We use Hashin's structure as a benchmark problem for our numerical computations; see \cite{dietrich_contrib:bib:H1962}. It consists of a coated circular inclusion in a matrix material, see Figure \ref{dietrich_contrib:fig:Hashin}.

Each material $\Omega_i$ is assumed to be isotropic. The Young's moduli and Poisson's ratios are denoted by $E_i$ and $\nu_i$, accordingly. In our tests, the radii are set to $r_1 = 0.25$ and $r_2 = 0.4$. If not mentioned otherwise, Young's moduli have the values $E_1 = \SI{100}{\giga\pascal}$ for the core material, $E_2 = \SI{1000}{\giga\pascal}$ for the coating and a resulting $E_3 = \SI{453.685}{\giga\pascal}$ for the matrix material, following the formulas found in \cite{dietrich_contrib:bib:KMS2015}.
Poisson's ratio is chosen to be $\nu = 0.3$ for all materials. A tolerance of $10^{-6}$ was used for the following computations.

\begin{figure}[bth]
\centering
\subfigure[$u_1(Y,y)=X_1(y)\vdotdot E(Y)$,\newline with $E(Y) = \text{Id}$]{\includegraphics[width=5.5cm]{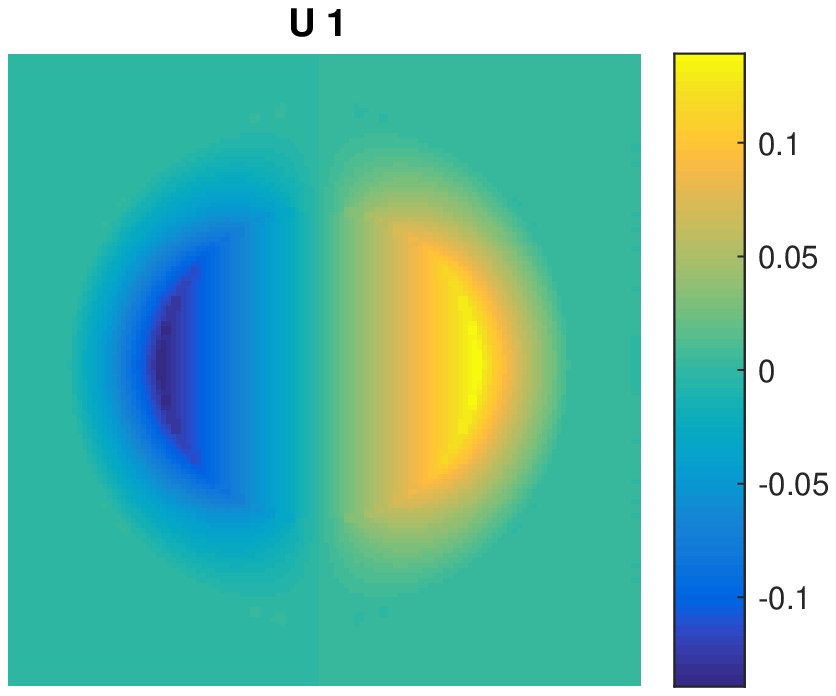}}
\hfill
\subfigure[$u_2(Y,y)=X_2(y)\righttherefore \nabla E(Y)$,\newline with $\left(\nabla E(Y)\right)_{111} = \left(\nabla E(Y)\right)_{222} = 1$,\newline $0$ otherwise]{\includegraphics[width=5.5cm]{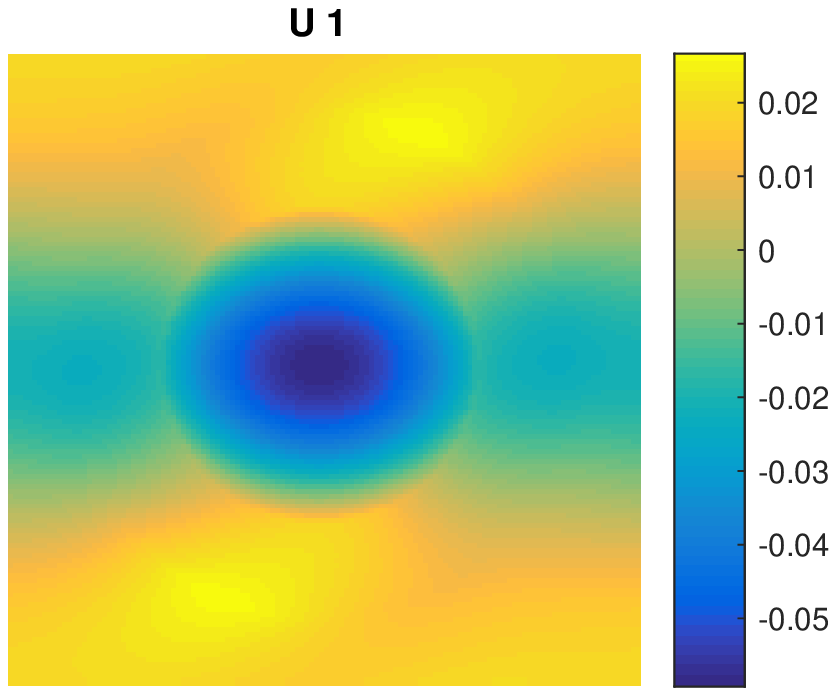}}
\caption{First component of displacement vectors for first (left) and second (right) order problems for the Hashin structure.}
\label{dietrich_contrib:fig:displacement_plots}
\end{figure}

In Figure \ref{dietrich_contrib:fig:displacement_plots}, the numerical solutions (using Algorithm \ref{dietrich_contrib:alg:acc_basic_scheme}) for the first component of the displacement vectors for the first and second order problems are shown. The underlying tensor grid consists of $128^2$ points. To solve the linear system in Algorithm \ref{dietrich_contrib:alg:acc_basic_scheme} we made use of \texttt{MATLAB}'s \texttt{bicgstab} function with a tolerance of $10^{-6}$.

\begin{figure}[bth]
\centering
\includegraphics[scale=1]{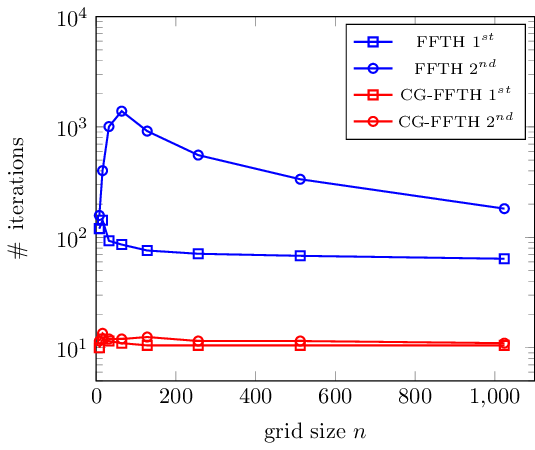}
\hfill
\includegraphics[scale=1]{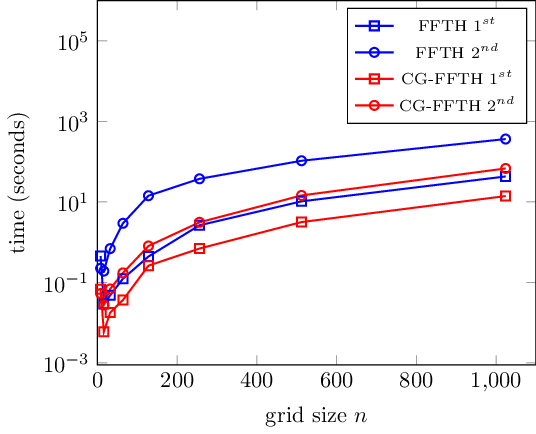}
\caption{Number of iterations and computation time needed on a grid of size $n\times n$.}
\label{dietrich_contrib:fig:versus_grid_size}
\end{figure}

Figure \ref{dietrich_contrib:fig:versus_grid_size} shows how the Basic Scheme --- shortened as \textit{FFTH} for \textit{FFT-based Homogenization} --- and its CG-version behave if the number of grid points gets larger. The standard algorithm needs much more iterations than the CG-version, especially for second order problems. It is important to note that the number of iterations is essentially independent of the grid size, although the Basic Scheme needs significantly more iterations for second order problems on smaller grids. The computation time for second order problems is in both algorithms noticeably higher than for first order problems. Due to the hierarchical structure of the problems, at least a factor of four was to be expected (three first order problems plus the second order problem itself). The computation of the polarization and body force terms result in additional overhead. A detailed comparison of the time ratios, i.e. the computation time of a second order problem divided by the time needed for the corresponding first order problem, can be found in Table \ref{dietrich_contrib:tab:time_ratio} for both algorithms. While the ratio keeps growing for the Basic Scheme, it appears to be limited for the CG-version around the expected value of four.

\begin{table}
\centering
\caption{Ratio of computation time for first and second order problem on a $(n\times n)$-grid.}
\label{dietrich_contrib:tab:time_ratio}
\begin{tabular}{p{2cm}p{2cm}p{2cm}}
\hline\noalign{\smallskip}
$n$ 	& FFTH 		& CG-FFTH  \\\midrule
$8$ 	& $0.4967$ 	& $0.7812$ \\
$16$ 	& $6.5105$ 	& $4.6337$ \\
$32$ 	& $14.4306$ 	& $3.8513$ \\
$64$ 	& $23.7463$ 	& $4.7257$ \\
$128$ 	& $32.3741$ 	& $3.0773$ \\
$256$ 	& $14.2606$ 	& $4.4829$ \\
$512$ 	& $10.1898$ 	& $4.5904$ \\
$1024$ 	& $8.4970$ 	& $4.8484$ \\
\noalign{\smallskip}\hline\noalign{\smallskip}
\end{tabular}
\end{table}

\begin{figure}[bth]
\centering
\includegraphics[scale=1]{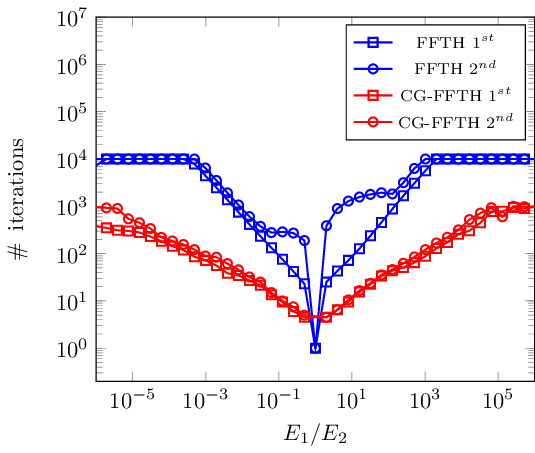}
\hfill
\includegraphics[scale=1]{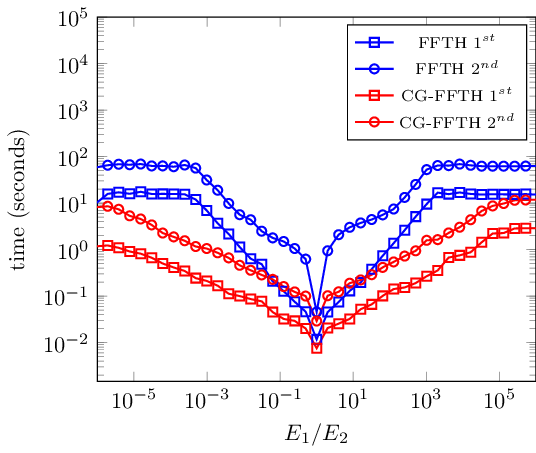}
\caption{Number of iterations and computation time plotted against the contrast of core and coating materials.}
\label{dietrich_contrib:fig:versus_contrast}
\end{figure}

For the plots shown in Figure \ref{dietrich_contrib:fig:versus_contrast}, we kept $E_1$ at a value of $\SI{100}{\giga\pascal}$ and changed the Young's modulus $E_2$ of the coating material. The computations were performed on a grid with $64^2$ points. In addition to the time gap between first and second order problems already shown before, we can see here that the number of iterations for the Basic Scheme surpasses $10^4$ iterations already for contrasts smaller than $10^{-3}$ or greater than $10^3$, whereas the CG-version can still handle these problems within a few hundred iterations. For the most part, its computation time is smaller as well.

\section{Conclusion}
\label{dietrich_contrib:sec:5}
Starting from an FFT-based scheme and its CG-version,
we have presented the generalization to higher-order derivatives and 
a comparison of the schemes for different orders in terms of number of iterations 
and computation time.
We are still working on an extensive quantitative analysis of the implications
on the effective properties, cf. \cite{dietrich_contrib:bib:B2007,dietrich_contrib:bib:MMS1999}.
This should be combined with a multiscale simulation using higher-order terms (FE-FFT coupling) \cite{dietrich_contrib:bib:KEtAl2017,dietrich_contrib:bib:SEtAl2014} in order
to obtain a meaningful assessment of the pros and cons of this approach.

\def\cprime{$'$}
\section*{Acknowledgments} 
The collaboration with H. Andr\"{a}, M. Kabel and M. Schneider, Fraunhofer ITWM Kaiserslautern, is gratefully acknowledged.

\ifx\undefined\bysame
\newcommand{\bysame}{\leavevmode\hbox to3em{\hrulefill}\,}
\fi

\end{document}